\documentclass[oneside,german,english,reqno]{amsart}
\usepackage{palatino}
\usepackage[T1]{fontenc}
\usepackage[latin1]{inputenc}
\usepackage{amssymb}
\IfFileExists{url.sty}{\usepackage{url}}
                      {\newcommand{\url}{\texttt}}

\makeatletter

\providecommand{\LyX}{L\kern-.1667em\lower.25em\hbox{Y}\kern-.125emX\@}

 \theoremstyle{plain}    
 \newtheorem{lem}{Lemma} 
 \theoremstyle{remark}
 \newtheorem*{rem*}{Remark}
 \theoremstyle{plain}    
 \newtheorem{thm}{Theorem} 

\DeclareMathOperator{\Vol}{Vol}
\DeclareMathOperator{\diam}{diam}
\DeclareMathOperator{\Haus}{Haus}
\DeclareMathOperator{\Boxx}{Box}
\DeclareMathOperator{\MST}{MST}

\usepackage{babel}
\makeatother
\begin{document}

\title{The minimal spanning tree and the upper box dimension}

\author{Gady Kozma}

\author{Zvi Lotker}

\author{Gideon Stupp}

\begin{abstract}
We show that the $\alpha $-weight of an MST over $n$ points in a
metric space with upper box dimension $d$ has a bound independent
of $n$ if $\alpha <d$ and does not have one if $\alpha >d$.
\end{abstract}
\maketitle

\section{Introduction}

The beginning point of this paper is the following old and well researched
theorem in Euclidean combinatorics

\begin{thm}
\label{thm:known}Let $v_{1},\dotsc ,v_{n}\in [0,1]^{d}$ and let
$T$ be the minimal spanning tree (MST) over $v_{1},\dotsc ,v_{n}$.
Then \[
E_{\alpha }(T):=\sum _{e\in T}||e||^{\alpha }\leq C(d)^{\alpha }n^{\max (0,1-\alpha /d)}.\]

\end{thm}
We note that the well known greedy algorithm for the minimal spanning
tree (see e.g.~\cite{CLR91}) shows that the spanning tree minimizing
$E_{1}(T)$ is also the spanning tree minimizing $E_{\alpha }(T)$
for any $\alpha $, and this fact is true in any graph, not just in
the Euclidean case.

This fact was a part of the general lore of the field for many years.
The case $\alpha =1$ was discussed in \cite{V51,FT40,F55,M84} with
increasingly improved constants. The case $\alpha =d$ (from which
the case of general $\alpha $ follows immediately, if one does not
care about optimizing the constant) was first proved in the case $d=2$
in \cite{GP68} and in the general case is usually attributed to \cite{SS87}
even though the result is not proved there explicitly. An almost trivial
argument using the existence of a space filling curve with precise
H\"older exponent can be found in the survey \cite[page 759]{S90}
and a different approach is presented in \cite{Y96,S01}. We cannot
resist citing the beautiful work done on the case of randomly distributed
points: see for example \cite{A96,ABNW99,AS}. However, we will skip
the vast literature on related problems of combinatoric optimization
such as the traveling salesperson, minimal matching, greedy matching,
Steiner minimal tree and so on.

Examine the proof of Gilbert and Pollak \cite{GP68}. Their approach
is as follows: place a rhombus, or diamond, on any edge of the MST
with an opening degree of $60^{\circ }$. It turns out that these
rhombi are disjoint which immediately proves the theorem with a constant
of $2\sqrt{3}$. The proof that they are disjoint uses specific facts
from the geometry of the plane. However, if one weakens the claim
to balls centered on the middle of each edge with radius $\frac{1}{10}$
the length of the edge, the proof becomes simpler and purely metric
(some assumption of uniform convexity is necessary). For the convinience
of the reader, we add a proof of this fact and hence of theorem \ref{thm:known}
in \S \ref{sub:classic} below.

Since we have a question phrased in metric terms, and a proof which
uses very little from the structure of $\mathbb{R}^{d}$, it seems
natural to try to extend the result to more general metric spaces,
and in particular to spaces with fractal dimensions. Thus we define
the {}``MST dimension'' using\begin{equation}
\dim _{\MST }(M):=\inf \{d:E_{d}(T)\leq C\quad \forall T\textrm{ a MST on points from }M\}.\label{eq:defdimMST}\end{equation}
We will prove that the MST dimension is equivalent to the upper box
dimension (also known as the upper Minkowsi dimension, or the entropy
dimension), namely

\begin{thm}
\label{thm:MSTBox}For every metric space $M$,\[
\dim _{\MST }(M)=\dim _{\Boxx }(M).\]

\end{thm}
See \cite{F90} for alternatives to the Hausdorff dimension and the
upper box dimension in particular. We note only that it is well known
that $\dim _{\Haus }\leq \dim _{\Boxx }$ and that they are equal
for a wide family of {}``well behaved'' sets.

One notable discrepancy remains between theorem \ref{thm:known} and
theorem \ref{thm:MSTBox} which is the critical case of $\alpha =d$.
Can theorem \ref{thm:MSTBox} be strengthened to show boundedness
in this case too? The short answer is no. The long answer is that
many factors come into play, most notably the connectivity of the
fractal. We plan to analyze specific cases in a separate paper.

\section{Proofs}

In the following, $C$ or $C(d)$ refer to constants, which could
be different from place to place.

\subsection{\label{sub:classic}The classic result}

\begin{lem}
\label{lem:geom}Let $w_{1},w_{2}\in \mathbb{R}^{d}$ satisfy $||w_{1}||\geq 1$,
$||w_{2}||\geq 1$ and $||w_{1}-w_{2}||\leq 1$. Then\[
\left\Vert \frac{w_{1}+w_{2}}{2}\right\Vert \geq \frac{\sqrt{3}}{2}.\]

\end{lem}
This is an easy exercise in plane geometry and we leave it to the
reader.

\begin{lem}
\label{lem:balls}Let $v_{1},\dotsc ,v_{n}\in [0,1]^{d}$ and let
$T$ be a spanning tree with minimal energy. For every edge $e\in T$,
let $B_{e}$ be a ball of radius $\frac{1}{10}||e||$ around the center
of $e$. Then the family of balls $\left\{ B_{e}\right\} _{e\in T}$
is pairwise disjoint.
\end{lem}
\begin{proof}
Assume to the contrary that $B_{e_{1}}\cap B_{e_{2}}\neq \emptyset $.
We may assume $||e_{1}||\geq ||e_{2}||$. Denote $e_{1}=(v_{1},v_{2})$
and $e_{2}=(w_{1},w_{2}$). Using the triangle inequality twice we
get \[
d\left(v_{1},\frac{w_{1}+w_{2}}{2}\right)\leq d\left(v_{1},\frac{v_{1}+v_{2}}{2}\right)+\frac{1}{10}||e_{1}||+\frac{1}{10}||e_{2}||<\frac{\sqrt{3}}{2}||e_{1}||.\]
We now use lemma \ref{lem:geom}. Since the conclusion of the lemma
fails, we know that one of it's assumptions fails. Since $||w_{1}-w_{2}||=||e_{2}||\leq ||e_{1}||$,
the failing assumption must be one of the first two, which means that
for some $i\in \{1,2\}$, $d(v_{1},w_{i})<||e_{1}||$. Similarly,
for some $j$, $d(v_{2},w_{j})<||e_{1}||$. Let $T'=T\setminus e_{1}$.
The removal of one edge splits the tree into two connected components.
Assume without loss of generality that $v_{1}$ and $w_{1}$ are in
the same component, and then $w_{2}$ would be in the same component.
Then $T''=T'\cup (v_{2},w_{j})$ is a spanning tree with lower energy
then $T$. This is a contradiction, and the lemma is proved.
\end{proof}

\begin{proof}
[Proof of theorem \ref{thm:known}]If $\alpha =d$ then this follows
immediately from lemma \ref{lem:balls}, since \begin{equation}
\sum _{e\in T}||e||^{d}\leq C(d)^{d}\sum _{e\in T}\Vol (B_{e})=C^{d}\Vol \left(\bigcup _{e\in T}B_{e}\right)\leq C^{d}\Vol ([-1,2]^{d})=C^{d}.\label{eq:d}\end{equation}
This also gives the required result for $\alpha >d$ since in this
case \begin{equation}
E_{\alpha }\leq \max ||e||^{\alpha -d}E_{d}\leq C(d)^{\alpha }.\label{eq:bigd}\end{equation}
 If $\alpha <d$ we can use H\"older's inequality (or $l^{p}-l^{q}$
duality) for $p=d/\alpha $ to get \begin{align*}
\sum ||e||^{\alpha } & =\sum e^{\alpha }\cdot 1\leq \Big |\Big |||e||^{\alpha }\Big |\Big |_{d/\alpha }\Big |\Big |\vec{1}\Big |\Big |_{d/(d-\alpha )}=\\
 & =\left(\sum ||e||^{d}\right)^{\alpha /d}n^{1-\alpha /d}\leq C(d)^{\alpha }n^{1-\alpha /d}.\qedhere 
\end{align*}

\end{proof}
\begin{rem*}
The dependency of $C$ on $d$ can be made explicit: it is enough
to take $C(d)=C\sqrt{d}$ where the right hand $C$ is absolute. This
value appears twice: in (\ref{eq:d}), where we need $C(d)>C\Vol (B^{d})^{-1/d}$,
and the estimate follows from the well known formula for the volume
of the $d$-dimensional ball. Next it appears in (\ref{eq:bigd}),
where we need $C(d)>C\diam [0,1]^{d}$ and this is clearly $>C\sqrt{d}$.
Thus we see that $C(d)=C\sqrt{d}$ is enough in both cases.
\end{rem*}

\subsection{The upper box dimension}

We shall use the following definition of the upper box dimension:
let $M$ be a metric space, and let $N(\epsilon )$ be the maximal
number of disjoint balls of radius $\epsilon $ in $M$. Then\[
\dim _{\Boxx }(M):=\overline{\lim _{\epsilon \rightarrow 0}}\frac{\log N(\epsilon )}{\log 1/\epsilon }.\]

\begin{lem}
\label{lem:MSTlarge}For every metric space $M$,\[
\dim _{\Boxx }(M)\leq \dim _{\MST }(M).\]

\end{lem}
\begin{proof}
If $\dim _{\Boxx }(M)=0$ there is nothing to prove. Therefore assume
$\alpha <\dim _{\Boxx }(M)$. By the definition of the box dimension,
it is possible to find sequences of disjoint balls $B(v_{1},\epsilon ),B(v_{2},\epsilon ),\dotsc ,B(v_{N(\epsilon )},\epsilon )$
such that \[
\overline{\lim _{\epsilon \rightarrow 0}}\epsilon ^{\alpha }N(\epsilon )=\infty .\]
Let $B(v_{i},\epsilon )$ be one such cover. Now examine the MST $T$
of the points $v_{i}$. We get that every edge of $T$ is $\geq \epsilon $
therefore $E_{\alpha }(T)\geq \epsilon ^{\alpha }N(\epsilon )$ and
since this is unbounded we get $\dim _{\MST }(M)\geq \alpha $. Since
this is true for every $\alpha <\dim _{\Boxx }(M)$ the lemma is proved.
\end{proof}
\begin{lem}
\label{lem:MSTsmall}For every metric space $M$,\[
\dim _{\MST }(M)\leq \dim _{\Boxx }(M).\]

\end{lem}
\begin{proof}
We may assume $M$ is bounded, since otherwise both dimensions are
$\infty $, and by scaling we may assume $\diam M=1$. Assume $\alpha >\dim _{\Boxx }(M)$
so that we have\begin{equation}
N(\epsilon )\leq C(1/\epsilon )^{\alpha }\quad \forall \epsilon >0\label{eq:Neps}\end{equation}
with $N(\epsilon )$ as above.

Recall the greedy algorithm for obtaining the MST: define $T_{0}=\{v_{0}\}$
and $T_{i}=T_{i-1}\cup \{v_{i}'\}$ where $v_{i}'$ is the vertex
not in $T_{i}$ closest to $T_{i}$, and the connection is via the
shortest edge connecting $v_{i}'$ to $T$. For every edge $e\in T$
satisfying $||e||>\epsilon $ define $B(e)$ to be a ball of radius
$\frac{1}{3}\epsilon $ centered at the vertex of $e$ which entered
last into $T$, or in other words, if $e=(v_{i}',v_{j}')$ and $i<j$
then $B(e)=B(v_{j}',\frac{1}{3}\epsilon )$. It follows that these
balls are disjoint: assume to the contrary that $B(e)\cap B(f)\neq \emptyset $,
and assume w.l.o.g that $e$ was added to the tree before $f$. Denote
$v$ the center of $B(e)$ and $w$ the center of $B(f)$ and then
$||v-w||<\frac{2}{3}\epsilon $. Then, when adding $f$, we have that
the partially constructed tree contains $v$, and we add $w$ which
contradicts the assumption $||f||>\epsilon $.

We now use (\ref{eq:Neps}), and get that there are no more than $C\epsilon ^{\alpha }$
edges of length $>\epsilon $ in $T$. Let $\alpha <\beta $. Then\begin{align*}
E_{\beta }(T) & =\sum _{e\in T}||e||^{\beta }=\sum _{k=0}^{\infty }\sum _{\substack{ e\in T\\
 2^{-k-1}<||e||\leq 2^{-k}}
}||e||^{\beta }\leq \\
 & \leq \sum _{k=0}^{\infty }C\left(2^{k+1}\right)^{\alpha }\left(2^{-k}\right)^{\beta }=C\sum _{k=0}^{\infty }2^{k(\alpha -\beta )}=C.
\end{align*}
 This holds for any $\beta >\dim _{\Boxx }(M)$ since we can always
find an $\alpha $ between $\dim _{\Boxx }(M)$ and $\beta $, and
the lemma is proved.
\end{proof}
\theoremstyle{remark}
\newtheorem*{remarks}{Remarks}
\begin{remarks}

\begin{enumerate}
\item The similarity between the proofs of lemma \ref{lem:MSTsmall} and
\ref{lem:balls} is evident. The only obstacle for us to use lemma
\ref{lem:balls} directly is that in a general space $M$ there is
no notion of a {}``point on the middle of an edge''. Even if we
embed $M$ into $\mathbb{R}^{d}$, the middle of the edge might not
belong to $M$, so that we cannot use it for the definition of the
upper box dimension.
\item In effect, we haven't used the greedy algorithm at all. What we have
shown is that the tree generated by the greedy algorithm has bounded
$\beta $-weight. Even if we didn't know that this tree is the MST,
it would definitely follow for the MST which has the minimal $\beta $-weight.
\item Theorem \ref{thm:MSTBox} may be generalized to spaces that satisfy
a weak triangle inequality, namely \begin{equation}
d(x,y)\leq C(d(x,z)+d(z,y)).\label{eq:weaktriangle}\end{equation}
The proof carries out without any significant change. We remark that
all the clasic equivalent definitions for the upper box dimension
(using covering or packing numbers), are also equivalent when (\ref{eq:weaktriangle})
replaces the triangle inequality.
\end{enumerate}
\end{remarks}

\end{document}